\documentclass[a4paper,11pt]{article}

\usepackage[margin=1in]{geometry}  

\usepackage{moreverb}

\usepackage[dvips,colorlinks,bookmarksopen,bookmarksnumbered,citecolor=red,urlcolor=red]{hyperref}

\usepackage{amsfonts,amsmath,amssymb,amsthm}
\usepackage{graphicx}

\usepackage{algorithm}
\usepackage{algorithmic}

\begin{document}

\title{RBF-PU Interpolation with Variable Subdomain Sizes and Shape Parameters}

\author{Roberto Cavoretto, Alessandra De Rossi, Emma Perracchione\footnote{Department of Mathematics \lq\lq G. Peano\rq\rq, University of Torino, via Carlo Alberto 10, I--10123 Torino, Italy. E-mails: roberto.cavoretto@unito.it, alessandra.derossi@unito.it, emma.perracchione@unito.it}}

\date{}

\maketitle

\begin{abstract}
	In this paper, we deal with the challenging computational issue of interpolating large data sets, with eventually non-homogeneous densities. To such scope, the Radial Basis Function Partition of Unity (RBF-PU) method has been proved to be a reliable numerical tool. However, there are not available techniques enabling us to efficiently select the sizes of the local PU subdomains which, together with the  value of the RBF shape parameter, greatly influence the accuracy of the final fit. Thus here, by minimizing an \emph{a priori} error estimate, we propose a RBF-PU method by suitably selecting variable shape parameters and subdomain sizes. Numerical results and applications show performaces of the interpolation technique.
\end{abstract}





\section{INTRODUCTION}

The interpolation via PU method consists in decomposing the domain into several \emph{subdomains} or \emph{patches} which, except for particular cases \cite{Safdari}, are always supposed to be hyperspheres of a fixed radius \cite{Cavoretto15b,Cavoretto16a}. Such fixed size of the subdomains penalizes the PU interpolant, especially when points with highly varying distributions are considered. Indeed, in these cases, problems as lack of information and/or ill-conditioning can arise.

In this work, considering  hyperspherical subdomains and focusing on RBFs as local approximants, the aim is to develop a method which enables us to select suitable radii for the PU subdomains and safe shape parameters of the local basis functions. For this purpose, we compute subsequent error estimates depending on these two quantities by selecting the optimal couple of values used to solve the local problem and keeping fixed the original data set. The term optimal is here used with abuse of notation; in fact, only if the function is known, the error can be exactly evaluated and thus the optimal values can be found without any uncertainty. Otherwise, all the techniques based on error estimates give approximated optimal values. 

The error estimates are computed with the use of a modified Leave One Out Cross Validation (LOOCV) scheme \cite{Golberg,Golub79,Rippa}; see also \cite[Chapter 14]{Fasshauer15}. More precisely, since our problem depends on two quantities, i.e. subdomain size and shape parameter, for each patch we perform a bivariate LOOCV. So the resulting method turns out to be accurate, and the use of such a flexible approach makes it particularly meaningful in real life problems. This follows from the fact that the computational issue consisting in approximating large and irregular data sets is rather common in a wide variety of applications. In order to point out accuracy and robustness of the new method, we provide a few numerical experiments, also investigating an application to Earth's topography.

The paper is organized as follows. In Section 2, after giving several preliminaries about the PU method, we focus on the construction of the new RBF-PU interpolant. Then, in Section 3 we provide numerical experiments and applications. Finally, Section 4 deals with conclusions.

\section{THE RBF-PU INTERPOLANT} 
\label{sec:PU}

Given a set ${ \cal X}_N= \{ \textbf{x}_i \in \mathbb{R}^{M},i=1, \ldots, N \}$ of $N$ distinct \textsl{data points}, also called \textsl{data sites} or \textsl{nodes}, in a domain $ \Omega \subseteq \mathbb{R}^{M}$, and a corresponding set $ {\cal F}_N= \{ f_i = f(\textbf{x}_i)  ,i=1, \ldots, N \}$ of \textsl{data values} or \textsl{function values} obtained by possibly sampling a function  $f: \Omega \longrightarrow \mathbb{R}$, the standard scattered data interpolation problem consists in recovering the function $f$ \cite{Iske11}.

In case of large scattered data sets, the PU method turns out to be extremely effective and easy to implement in any dimension. Precisely, we consider a partition of the open and bounded domain $ \Omega$ into $d$ subdomains $ \Omega_j$,  such that $ \Omega  \subseteq \cup_{j=1}^{d} \Omega_j$,  with some mild overlap among them \cite[Chapter 15]{Wendland05}. Furthermore, we select a family of compactly supported, non-negative, continuous functions $W_j$, with $\textrm{supp}(W_j) \subseteq \Omega_j$ and forming a partition of unity, i.e. $\sum_{j=1}^{d} W_j(\textbf{x}) = 1,$ $\textbf{x} \in \Omega$. Then, the global interpolant ${\cal I}$ is formed by the weighted sum of $d$ local RBF approximants $R_j$ \cite{Cavoretto16b}, i.e.
\begin{eqnarray} \label{intg}
{\cal I}\left( \textbf{x}\right)= \sum_{j=1}^{d} R_j\left( \textbf{x} \right) W_j \left( \textbf{x}\right), \quad \textbf{x} \in \Omega, 
\end{eqnarray}
with
\begin{eqnarray} \label{defrbf} 
R_j( \textbf{x} )= \sum_{k=1}^{N_j} c_k^j \phi_{\varepsilon_j} ( || \textbf{x} -  \textbf{x}^j_k ||_2 ), 
\end{eqnarray}
where $\phi_{\varepsilon_j} :[0,\infty) \rightarrow \mathbb{R}$ is a RBF of shape parameter $\varepsilon_j$,  $N_j$ indicates the number of data points belonging to $\Omega_j$ and $\textbf{x}_k^j \in {\cal X}_{N_j}= {\cal X}_N \cap \Omega_j$, with $k=1, \ldots, N_j$.

The coefficients $\{ c_k^j \}_{k=1}^{N_j}$ are determined by imposing $R_j (\textbf{x}_i^j )=f_i^j,$ $i=1,\ldots,N_j.$
Thus, in order to find the PU interpolant ${\cal I}$, we need to solve $d$ linear systems of the form
\begin{eqnarray*}
A_j \textbf{c}_j= \textbf{f}_j ,
\end{eqnarray*}
where $(A_j)_{ik}= \phi_{\varepsilon_j}  ( || \textbf{x}^j_i - \textbf{x}_k^j ||_2 ),$  $i,k=1, \ldots, N_j, $
$  \textbf{c}_j=  (c_1^j, \ldots, c_{N_j}^j  )^T$ and $  \textbf{f}_j = (f_1^j, \ldots , f_{N_j}^j  )^T$.

As evident from (\ref{defrbf}), the accuracy of the PU fit depends on both $N_j$, i.e. the radius of the $j$-th patch $\delta_j$, and the shape parameter $\varepsilon_j$. Usually, they are supposed to be fixed for all $\Omega_j$ (see e.g. \cite[Chapter 29]{Fasshauer}), i.e. $\varepsilon = \varepsilon_j$ and $\delta = \delta_j$, where
\begin{eqnarray}
\delta_j = \frac{1}{d^{1/M}}, \quad  j=1, \ldots ,d.
\label{fisso}
\end{eqnarray} 

Here instead, $\delta_j$ and $\varepsilon_j$ are supposed to vary among the subdomains and selected by means of  cross validation schemes \cite[Chapter 17]{Fasshauer}, properly modified for bivariate optimization problems. Specifically, we focus on the LOOCV, firstly introduced in \cite{Golub79} and  further developed in \cite{Rippa}. Such approach is always performed to find the optimal value of the shape parameter for a global interpolation problem, while here we are interested in selecting  the optimal couple $(\delta_j,\varepsilon_j)$. 

Let us consider an interpolation problem on $\Omega_j$ of the form (\ref{intg}) and, for a fixed $i \in \{1, \ldots, N_j\}$, let
\begin{eqnarray*}
R^{(i)}_j  ( \textbf{x} )=  \sum_{k=1, k \neq i}^{N_j} c_k^j \phi_{\varepsilon_j}  ( || \textbf{x} -\textbf{x}^j_k ||_2 ) 
\end{eqnarray*}
and
\begin{eqnarray*}
e^j_i= f^j_i-R^{(i)}_j (\textbf{x}^j_i)
\end{eqnarray*}
be respectively the $j$-th interpolant obtained leaving out the $i$-th data on $\Omega_j$ and the  error at the $i$-th point. Then, following \cite{Rippa} and \cite[Chapter 17]{Fasshauer}, in order to obtain an error estimate,  we compute 
\begin{eqnarray}
\textbf{e}_j = \left( e^j_1, \ldots, e^j_{N_j} \right)= \left( \frac{c^j_1}{ \left(A_j \right)_{11}^{-1}}, \ldots, \frac{c^j_{N_j }}{\left(A_j \right)_{N_j N_j}^{-1}} \right),
\label{errori}
\end{eqnarray}
where $c^j_i$ is the $i$-th coefficient of the local RBF interpolant $R_j$ based on the full data set and $(A_j )_{ii}^{-1}$ is the $i$-th diagonal element of the inverse of the corresponding local interpolation matrix. 

In order to select the optimal couple $(\delta_j, \varepsilon_j)$ for each PU subdomain,  we compute (\ref{errori}) for  several values of the radius  $(\delta_{j_1}, \ldots, \delta_{j_P})$ and of the shape  parameter $(\varepsilon_{j_1}, \ldots, \varepsilon_{j_Q})$.

In what follows, to stress the dependence of (\ref{errori}) also from the shape parameter, for a fixed  $p \in \{1, \ldots, P \} $ and a fixed $q \in \{1, \ldots, Q \} $, we use the notation 
\begin{eqnarray*}
\textbf{e}_j \left(\delta_{j_p},\varepsilon_{j_q} \right) = \left( e_1^j \left( \delta_{j_p},\varepsilon_{j_q} \right), \ldots, e^j_{N_j} \left( \delta_{j_p},\varepsilon_{j_q} \right) \right).
\end{eqnarray*}
Without any restrictions we focus on the maximum norm and we define
\begin{equation}
E_j=
\left(
\begin{array}{ccc}
||\textbf{e}_j (\delta_{j_1},\varepsilon_{j_1}) ||_{\infty} & \cdots & ||\textbf{e}_j (\delta_{j_1},\varepsilon_{j_Q}) ||_{\infty} \\	
\vdots & \ddots & \vdots \\
||\textbf{e}_j (\delta_{j_P},\varepsilon_{j_1}) ||_{\infty} & \cdots & ||\textbf{e}_j (\delta_{j_P},\varepsilon_{j_Q}) ||_{\infty} \\	
\end{array}
\right). 
\label{mate}
\end{equation}
The elements stored in the matrix (\ref{mate}) are error estimates for several values of the PU radius and of the RBF shape parameter. Therefore, we compute the $j$-th local approximant with the couple of values $(\delta_j, \varepsilon_j)$ such that
\begin{eqnarray*}
||\textbf{e}_j (\delta_{j},\varepsilon_{j}) ||_{\infty} = \min_{ p=1, \ldots, P} \left( \min_{  q=1, \ldots, Q} (E_j)_{pq} \right).
\label{opt_coppia}
\end{eqnarray*}

Trivially, we can observe that the searching intervals for both the radius and the shape parameter need to be properly selected. Since many researchers have already worked on the problem of finding suitable values for shape parameter, we can easily guess how to select a good range for it, see for instance \cite{Driscoll-Fornberg02,FasshauerZhang,Fornberg11}. In other words, for what concerns the shape parameter the notation simplifies, since for each subdomain we can consider the same interval, namely  $(\varepsilon_{1}, \ldots, \varepsilon_{Q})$.
While, for the PU radius $\delta_j$, we search its optimal value in an interval of the form
\begin{equation}
(\delta_{j_1}, \ldots, h \delta_{j_1}), \quad h \in \mathbb{R}^{+},  \quad h>1, \quad \textrm{and} \quad  \delta_{j_1} \quad \textrm{is such that} \quad Card(\Omega_j)  \geq N B(\delta_{j_1}),
\label{range}
\end{equation}
where $B(\delta_{j_1})$ is the hypervolume of the hypersphere of radius $\delta_{j_1}$. In this way, we avoid the problem of considering empty  subdomains.
Moreover, since the upper bound of the searching interval (\ref{range}) is proportional to the lower bound and since this lower bound is large only if the density of points is low, we also avoid problems arising from high density of points, i.e. systems are not \emph{too large} and the ill-conditioning is kept under control. Finally, we can note that, aside from the complexity cost of the RBF-PU method with fixed radii and shape parameters \cite{Cavoretto15b}, we also require the computation of the matrix inverse for each $\delta_{j_p}$ and $\varepsilon_{j_q}$.

\section{NUMERICAL EXPERIMENTS AND APPLICATIONS}
\label{sec:ne}
This section is devoted to show, by means of  numerical simulations, the flexibility and the accuracy of the proposed RBF-PU method. Tests are carried  out considering  the so-called \emph{product} function
\begin{eqnarray*}
f(x_1,x_2)  = 16 x_1 x_2 (1-x_1) (1-x_2).
\end{eqnarray*}
As interpolation points, we take uniformly random Halton data on $\Omega=[0,1]^2$.

To point out the accuracy of the new interpolant, we compute the Maximum Absolute Error (MAE) and the Root Mean Square Error (RMSE):
\begin{eqnarray*} \label{MAE}
\textrm{MAE} = \max_{ i=1, \ldots, s} |f(\tilde{\textbf{x}}_i) - {\cal I}(\tilde{\textbf{x}}_i) |, \qquad \textrm{RMSE} = \sqrt{\frac{1}{s}\sum_{i=1}^{s} |f(\tilde{\textbf{x}}_i) - {\cal I}(\tilde{\textbf{x}}_i)|^2},
\end{eqnarray*}
where $\tilde{\textbf{x}}_i$, $i=1, \ldots, s$, is a grid of $40 \times 40$ points in which the interpolant is sought. 
As basis functions we consider the Mat\'ern $C^2$ and  Inverse MultiQuadric (IMQ) functions \cite[Chapter 4]{Fasshauer}. They are respectively defined as:
\begin{eqnarray*}
\phi_{\varepsilon}^1  (r) = e^{- \varepsilon r} (1+\varepsilon r) \quad \textrm{and} \quad \phi_{\varepsilon}^2 (r)=\left(1+(\varepsilon r)^2 \right)^{-1/2},
\label{IMQ}
\end{eqnarray*}
where $r$ is the Euclidean norm.

In Table \ref{tabeu1} we show the results obtained by means of the RBF-PU interpolant with variable values of $\delta_j$ and $\varepsilon_j$. Specifically, we take the IMQ and we choose $30$ values for the shape parameter in the range $(0.001, 10)$.  Moreover, we fix the initial ranges for the radii as in (\ref{range}), with $h=2$ and $P=6$. As a comparison, we also report  the errors of the RBF-PU method obtained by fixing the radius as in (\ref{fisso}) and the shape parameter $\varepsilon=0.6$.

\begin{table}[ht!]
\begin{center}
		\begin{tabular}{ccccc} 	\hline\noalign{\smallskip}
			$N$ 	  & RMSE$_{(\delta_j,\varepsilon_j)}$  & MAE$_{(\delta_j,\varepsilon_j)}$ & RMSE$_{(\delta,\varepsilon)}$  & MAE$_{(\delta,\varepsilon)}$\\
			\hline 
 			\rule[0mm]{0mm}{3ex}
 			$\hskip-2pt 289$  &  $1.03{\rm E}-05$  &  $2.36{\rm E}-04$  &  $3.64{\rm E}-03$  &  $5.66{\rm E}-02$ 	   \\
 			\rule[0mm]{0mm}{3ex}
 			$1089$      & $2.88{\rm E}-06$  &  $7.89{\rm E}-05$	& $7.57{\rm E}-04$  &  $1.52{\rm E}-02$     \\
 			\rule[0mm]{0mm}{3ex}
 			$4225$      & $3.84{\rm E}-07$  &  $1.39{\rm E}-05$ & $3.88{\rm E}-04$  &  $1.01{\rm E}-02$ 	  \\
 			\rule[0mm]{0mm}{3ex}
 			$16641$     & $9.67{\rm E}-08$  &  $3.15{\rm E}-06$ & $8.27{\rm E}-04$  &  $3.27{\rm E}-02$     \\
 			\rule[0mm]{0mm}{3ex}
 			$66049$     & $2.68{\rm E}-08$  &  $6.80{\rm E}-07$ & $1.08{\rm E}-05$  &  $1.09{\rm E}-04$     \\
 			\hline 
		\end{tabular}
		\end{center}
	\caption{RMSEs and MAEs computed on Halton points via the RBF-PU methods by using the IMQ as local RBF interpolant.}
	\label{tabeu1}
\end{table}


Finally, we test the RBF-PU scheme with real world data. In particular, we consider the so-called \emph{glacier} data set. It consists of $8345$ points (with non-homogeneous density) representing digitized height contours of a glacier \cite{Wendland01}. The difference between the highest and the lowest point  is $800$ m. Because of the high variability of the points, we use as local approximant the Mat\'ern $C^2$. To test accuracy of the variable RBF-PU method, since in this practical situation we cannot sample data from a function, we randomly take $90$ points of the glacier data set as validation points, thus obtaining $\textrm{RMSE}=0.65$ m and $\textrm{MAE}=3.31$ m. 


\section{CONCLUDING REMARKS}
\label{sec:concl}
The proposed PU interpolant is based on safely selecting the parameters affecting the accuracy of the local fits. Numerical evidence and a brief sketch of an Earth's topography application show that the proposed method turns out to be accurate also when irregular data are considered. 

Work in progress consists in considering further possible shapes for the PU subdomains. This is not trivial since several requirements for the covering might be not easily satisfied. 

\section{ACKNOWLEDGMENTS}
The authors sincerely thank the two anonymous referees for helping to improve the paper. This work was partially supported by the University of Torino via grant \lq\lq Metodi numerici nelle scienze applicate\rq\rq, 2014.


\nocite{*}
\bibliographystyle{elsart-num-sort}%
\bibliography{bib_cavoretto}%

\begin{thebibliography}{10}
\expandafter\ifx\csname url\endcsname\relax
  \def\url#1{\texttt{#1}}\fi
\expandafter\ifx\csname urlprefix\endcsname\relax\def\urlprefix{URL }\fi

\bibitem{Fornberg11}
{B. Fornberg}, {E. Larsson}, {N. Flyer}, SIAM J. Sci. Comput. 33 (2011)
  869--892.

\bibitem{Fasshauer}
{G.E. Fasshauer}, Meshfree approximation methods with \textsc{Matlab}, World
  Scientific, Singapore, 2007.

\bibitem{FasshauerZhang}
{G.E. Fasshauer}, {J.G. Zhang}, Numer. Algorithms 45 (2007) 345--368.

\bibitem{Fasshauer15}
{G.E. Fasshauer}, {M.J. McCourt}, Kernel-based approximation methods using
  \textsc{Matlab}, World Scientific, Singapore, 2015.

\bibitem{Golub79}
{G.H. Golub}, {M. Heath}, {G. Wahba}, Technometrics 21 (1979) 215--223.

\bibitem{Iske11}
A.~Iske, Rend. Sem. Mat. Univ. Pol. Torino 69 (2011) 217--246.

\bibitem{Golberg}
{M.A. Golberg}, {C.S. Chen}, {S.R. Karur}, Eng. Anal. Bound. Elem. 18 (1996)
  9--17.

\bibitem{Cavoretto15b}
{R. Cavoretto}, {A. De Rossi}, SIAM J. Sci. Comput. 37 (2015) A1891--A1908.

\bibitem{Cavoretto16a}
{R. Cavoretto}, {A. De Rossi}, {E. Perracchione}, J. Sci. Comput. 68 (2016)
  395--415.

\bibitem{Cavoretto16b}
{R. Cavoretto}, {A. De Rossi}, {E. Perracchione}, Comput. Math. Appl. 71 (2016)
  2568--2584.

\bibitem{Rippa}
S.~Rippa, Adv. Comput. Math. 11 (1999) 193--210.

\bibitem{Safdari}
A.~Safdari-Vaighani, A.~Heryudono, E.~Larsson, J. Sci. Comput. 64 (2015)
  341--367.

\bibitem{Driscoll-Fornberg02}
{T.A. Driscoll}, {B. Fornberg}, Comput. Math. Appl. 43 (2002) 413--422.

\bibitem{Wendland01}
H.~Wendland, IMA J. Numer. Anal. 21 (2001) 285--300.

\bibitem{Wendland05}
H.~Wendland, Scattered data approximation, Cambridge Univ. Press, Cambridge,
  2005.

\end{thebibliography}

\end{document}